\documentclass{amsart}

\usepackage{amsmath}
\usepackage{verbatim}
\usepackage{upref}
\usepackage{amsfonts,amssymb}

\newtheorem{theorem}{Theorem}[section]
\newtheorem{lemma}[theorem]{Lemma}
\newtheorem{proposition}[theorem]{Proposition}
\newtheorem{corollary}[theorem]{Corollary}

\newtheorem{main lemma}[theorem]{Main Lemma}

\theoremstyle{definition}
\newtheorem{definition}[theorem]{Definition}

\newtheorem{remark}[theorem]{Remark}

\newtheorem{notation}[theorem]{Notation}

\newcommand{\N}{\mathbb{N}}
\newcommand{\C}{\mathbb{C}}
\newcommand{\Z}{\mathbb{Z}}
\newcommand{\T}{\mathbb{T}}
\newcommand{\A}{C^\circ}

\newcommand{\B}{\mathbb{B}}

\newcommand{\id}{\text{\rm id}}
\newcommand{\Alg}{\text{\rm Alg}}

\newcommand{\alg}{\text{\rm Alg$(H,k)$}}
\newcommand{\Hom}{\text{\rm Hom}}
\newcommand{\End}{\text{\rm End}}
\newcommand{\Super}{\text{\rm Super}}

\newcommand{\rp}{\text{\rm Rep}}
\newcommand{\ev}{\text{\rm ev}}
\newcommand{\aut}{\text{\rm Aut}}
\newcommand{\im}{\text{\rm Im}}

\newcommand{\nnn}{\overline{\N}}

\newcommand{\w}{\text{\rm w}}

\begin{document}

\title[Discrete semilattices grading algebras]
{Action of Pontryagin dual of semilattices grading algebras}

\author[Lucio Centrone]{Lucio Centrone$^*$}

\address{$^*$ IMECC, Universidade Estadual de Campinas (UNICAMP), Rua S\'ergio Buarque de Holanda, 651, Cidade Universit\'aria ``Zeferino Vaz'', Distr. Bar\~ao geraldo, Campinas, S\~ao Paulo, Brazil.}
\email{centrone@ime.unicamp.br}




\begin{abstract}
Let $A$ be a unitary algebra and $G$ be a finite abelian group. Then a $G$-graded algebra is merely a $G$-algebra and viceversa because of the fact that $G$ and its group of characters $\widehat{G}$ are isomorphic. This fact is no longer true if we substitute $G$ with infinite or non-abelian groups. In this paper we try to obtain similar results for a special class of abelian monoids, i.e, the bounded semmilattices. Moreover, if $S$ is such a monoid, we are going to investigate the role of $S$ and its Pontryagin dual $\widehat{S}$ over the algebra $A$, in the case $A$ is $S$-graded.
\end{abstract}

\maketitle

\section{Introduction}
If $G$ is a group, we say that an algebra $A$ is $G$-graded if $A=\bigoplus_{g\in G}A^g$, where $A^g\subseteq A$ are subspaces, and for every $g,h\in G$, $A^gA^h\subseteq A^{gh}$.
It is well known that if $G$ is a finite abelian group and $A$ is a unitary $G$-graded
algebra, then $G$
acts on $A$. More precisely, we consider the group $\widehat{G}$ of characters
of $G$, then $\widehat{G}$ acts on $A$. In fact, we may set
$\chi\cdot a=\chi(g)a$, provided $a\in A^g$, and it turns out the previous
relation is a group action over $A$. Conversely, if $\widehat{G}$ acts on $A$,
$A$ is $G$-graded once we set for any $g\in G$,
$A^g= \{a\in A|\chi\cdot a=\chi(g)a,\forall\chi\in\widehat{G}\}$.
In the case $G$ is abelian, we have $G\cong\widehat{G}$, then $A$
is a $G$-graded algebra if and only if $G$ acts on $A$.
In \cite{com1}, the authors recognized that for finite groups $G$, 
a coaction of $G$ on a algebra $A$ is merely a $G$-grading of $A$.
We recall that a $G$-grading of an algebra $A$ over
a field $k$ is equivalent to having a $(kG)^*$-module algebra action
of $(kG)^*$, the classical dual space of the group algebra $kG$.

Following this line of research, in this paper we study
what happens for algebras graded by monoids instead of groups.
Let $S$ be a monoid, then starting with an $S$-graded algebra
$A$, is it still true that there exists a monoid action over $A$?
Is there any relation between the monoid $S$ and its Pontryagin dual $\widehat{S}$? In this work we try to answer, at least partially, to both of the questions. The paper is divided
into three parts. The first one has two aims: recalling the general notions
of coalgebra, bialgebra and their related structures and reproducing
the duality theorems of \cite{com1}. In particular, we prove that if
$A$ is an algebra graded by a monoid $S$, there exists an action of
the group-like elements of the finite dual algebra
of the monoid algebra $kS$ (see Theorem \ref{gradingaction} and Corollary \ref{classic1}). In the second section, we recall some
topological definitions and tools in order to give a better description
of the monoids acting on a $S$-graded algebra. In few words, we focus on
the fact that the group-like elements of the finite dual of $kS$ and the
Pontryagin dual of $S$ are not so ``far'' in the case $S$ is a discrete bounded semilattice (see Theorem \ref{repmonoid} and Corollary \ref{classic2}). In the tird section we apply the results of the
previous two sections on concrete algebras. We consider a finite semilattice grading the upper triangular matrices. Moreover, we study
the case of the algebra $P$ of letterplaces in a countable set of
indeterminates. We consider the monoid $\nnn=\N\cup\{-\infty\}$, endowed
with the $\max$ binary operation and a special $\nnn$-grading on $P$.
It turns out that $\underline{\N}=\nnn\cup\{+\infty\}$, endowed with
the $\min$ binary operation, acts on $P$ and we describe explicitely
the finite dual and the group-like elements of the monoid algebra $k\nnn$. It turns out that in the first example we obtain a duality result similar to that of group grading algebras while in the second example the dualism fails.

\section{Definitions and preliminary results}

\subsection{Coalgebras and group-like elements}

Assume $k$ is a field of characteristic zero. Notice that the previous hypothesis is not necessary in many of the following facts. Indeed we will underline when the characteristic 0 becomes crucial. 

\begin{notation}
Let $(A,m,u)$ be a unital associative $k$-algebra, where $m:A\otimes A\to A$
is the multiplication map and $u:k\to A$ is the unital map. We denote by
$\End_k(A)$ the algebra of $k$-linear endomorphisms of $A$ (as a
vector space) and by $\aut_k(A)$ the group of $k$-linear automorphisms
of $A$. Moreover, we define $\End(A)$, $\aut(A)$ respectively the
monoid (group) of algebra endomorphisms (automorphisms) of $A$.
Let $(S,\cdot)$ be a monoid with identity element $e$.
We denote by the symbol $k S$ the monoid algebra generated by $S$.

Notice that we use the symbol $e$ for the unit element of a monoid $S$ and the symbol 1 for the usual unit of a $k$-algebra. By the way, in the monoid algebra $kS$ the unit 1 coincides with the element $1\cdot e$. Then we are going to use (with abuse of notation) either the symbol $e$ or the symbol 1 in order to denote the unit of $kS$. 
\end{notation}

\begin{definition}
Let $(S,\cdot)$ be a monoid. We say that $S$ acts on a $k$-algebra $A$
or that $A$ is an $S$-algebra if there is a monoid homomorphism
$\alpha:S\rightarrow \End(A)$. In what follows we denote by $\alpha_s$
the image of $s$ in $\End(A)$. In the same way, if $(G,\cdot)$
is a group then $A$ is a $G$-algebra if one has a group homomorphism
$\alpha:G\rightarrow \aut(A)$.
\end{definition}

We recall the main definitions and results about coalgebras and bialgebras and its related structure. We shall indicate by $(C,\Delta,\epsilon)$ the $k$-coalgebra $C$ equipped with the coalgebra maps $\Delta$ and $\epsilon$. Indeed, we shall use the same notation for any $k$-bialgebra, too. A subspace $I$ of $C$ is called coideal if $\Delta(I)\subseteq I\otimes C +
C\otimes I$ and $\epsilon(I) = 0$. In this case, the quotient space $C/I$
becomes a coalgebra in a natural way. If $(C,\Delta_C,\epsilon_C)$ and $(D,\Delta_D,\epsilon_D)$ are two
$k$-coalgebras, we say that a map $f:C\rightarrow D$ is a coalgebra morphism if $\Delta_D\circ f =
(f\otimes f)\circ\Delta_C$ and $\epsilon_C = \epsilon_D\circ f$.
Clearly, the kernel of every coalgebra morphism $f:C\rightarrow D$ is a coideal in $C$
and the image is a subcoalgebra of $D$. Moreover, one has that $C/\ker(f)$ and
$\im(f)$ are isomorphic coalgebras.

We consider the following definiton.
\begin{definition}
Let $(C,\Delta,\epsilon)$ be a $k$-coalgebra and $c\in C$. We say that $c$
is a group-like element if $\Delta(c) = c\otimes c$ and $\epsilon(c) = 1$.
We shall denote by $G(C)$ the set of group-like elements of $C$.
\end{definition}

It is well knwon if $(C,\Delta,\epsilon)$ be a $k$-coalgebra, the elements of $G(C)$
are linearly independent, then if $C=kS$, $G(C)=S$. In what follows, we shall need the following lemma.

\begin{lemma}
\label{grouplikequotient}
Let $S$ be a monoid and let $I$ be a coideal of $k S$. Then
$G(k S/I) = S + I$, where $S + I = \{ s+I | s\in S\}$.
\end{lemma}

\proof
It is clear that $S + I\subseteq G(k S/I) = G$. Consider now
$a + I\in G$, where $a = \sum\alpha_s s$.
By the definition of quotient coalgebra, we have that $\Delta'(a+I) =
\sum\alpha_s(s+I)\otimes(s+I)$. On the other hand,
$a+I\in G$ implies that $\Delta'(a+I)=(a+I)\otimes(a+I) =
\sum_{r,t\in S}\alpha_r\alpha_t(r+I)\otimes(t+I)$. Comparing the two expressions
for $\Delta'(a+I)$ one obtains that $\Delta'(a+I) =
\sum\alpha_s^2(s+I)\otimes(s+I)$, then $\alpha_s^2=\alpha_s$. The previous relation and the fact that $k$ is a field of characteristic 0 lead us to say that $\alpha_s=1$ and we are done.
\endproof

We recall that for any vector space $V$ over $k$, $V$ and its linear dual $V^*$ determine a non-degenerate bilinear form $\langle\ ,\ \rangle:V^*\otimes V\rightarrow k$ such that $\langle f,v\rangle=f(v)$ for any $f\in V^*$ and $v\in V$. Moreover, if $W$ is a vector space over $k$ and $\psi:V\rightarrow W$ is $k$-linear, we denote by $\psi^*$ the \textit{transpose} of $\psi$. More precisely, $\psi^*:W^*\rightarrow V^*$ such that \[\psi^*(f)(v)=f(\psi(v)),\] for all $f\in W^*$, $v\in V$.

We have the following (see \cite{mon1}, Lemma 1.2.2):

\begin{lemma}
If $(C,\Delta,\epsilon)$ is a $k$-coalgebra, then $C^*$ is a $k$-algebra with multiplication $m=\Delta^*$ and unit $u=\epsilon^*$.
\end{lemma}

We consider now the finite dual coalgebra of an algebra. 

\begin{definition}
\label{finitedual}
Let $(A,m,u)$ be a $k$-algebra and consider $A^* = \Hom_k(A,k)$, the linear dual
of $A$. Define \[A^\circ =\{f\in A^* \mid f(I)=0,\ \mbox{for some ideal}\
I\subseteq A\ \mbox{s.t.}\ \dim(A/I) < \infty\}.\] Then, one has that
$(A^\circ,\Delta,\epsilon)$ is a $k$-coalgebra (see \cite{mon1}
Proposition 9.1.2) where $\Delta f(a\otimes b) = m^* f(a\otimes b) = f(a b)$,
for all $f\in A^\circ$ and $a,b\in A$. Moreover, one defines $\epsilon =
u^*:A^\circ\to k$. The coalgebra $A^\circ$ is called the finite dual of
the algebra $A$.
\end{definition}

A key result is the following one (see \cite{mon1}, Example 9.1.4).

\begin{proposition}
\label{grouplikealgebra}
Let $A$ be a $k$-algebra and denote $\Alg(A,k)$ the monoid of algebra
homomorphisms $A\to k$. One has that $G(A^\circ) = \Alg(A,k)$.
\end{proposition}

If $C=kS$ is a $k$-coalgebra, the notion of graded vector space and the notion of $C$-modules coincide in a certain sense. We have the following results.

\begin{proposition}
\label{comodulegrading}
Let $(S,\cdot)$ be a monoid, then we consider the $k$-coalgebra $C = kS$. Then $M$ is a right $C$-comodule if and only if
$M$ is a $S$-graded vector space.
\end{proposition}

\begin{proposition}
\label{comodulemodule}
Let $A$ be a $k$-algebra and let $C$ be a $k$-coalgebra. If $M$ is a
vector space over $k$ then we have the following.
\begin{itemize}
\item[(i)] If $M$ is a right $C$-comodule then $M$ is also a left $C^*$-module.
\item[(ii)] Let $M$ be a left $A$-module. Then $M$ is a right $A^\circ$-comodule
if and only if for any $m\in M$ one has that $\dim\{A\cdot m\} < \infty$.
\end{itemize}
\end{proposition}

About (i) of the previous result, note that if $\rho:M\rightarrow M\otimes C$
is the comodule map such that $\rho(m) = \sum m_0\otimes m_1$ for any $m\in M$
and if $f\in C^*$ then $M$ becomes a left $C^*$-module under the map
$\gamma:C^*\otimes M\rightarrow M$ such that $\gamma(f\otimes m) =
\sum \langle f,m_1\rangle m_0$. Therefore, one obtains immediately
the following result.

\begin{corollary}
\label{modulestructure}
Let $(S,\cdot)$ be a monoid. Moreover, let $M$ be a $S$-graded vector space
and consider the monoid algebra $C = k S$ endowed with the structure
of $k$-coalgebra. Then $M$ is a left $C^*$-module under the map
$\gamma:C^*\otimes M\rightarrow M$ such that $\gamma(f\otimes m) =
\sum \langle f, s\rangle m_s$ for any $f\in C^*$ and $m = \sum_s m_s\in M$.
\end{corollary}

\subsection{Bialgebras and action of group-like elements}

We focus on bialgebras and how the group-like elements act on
a $k$-algebra graded by a monoid. Note that if $(B,\Delta,\epsilon)$ is a $k$-bialgebra then $G(B)$ is a monoid under
the multiplication of $B$.

The following result generalizes the action of a monoid
over an algebra.

\begin{proposition}
\label{actionmodulealgebra}
Let $S$ be a monoid and $A$ a $k$-algebra. Then $S$ acts on $A$ if and only if $A$
is a left $k S$-module algebra.
\end{proposition}

Now we state one classical result that we shall use later on. For the proof of the, see \cite{mon1} Theorem 9.1.3.
\begin{theorem}
\label{bialgebra}
Let $(B,\Delta,\epsilon)$ be a $k$-bialgebra. Then, the finite dual $B^\circ$
is also a $k$-bialgebra.
\end{theorem}

In particular, let us consider the monoid algebra $C=k S$, its linear dual $C^*$
and its finite dual $\A\subseteq C^*$. We have that $\A$ is an algebra with
the usual pointwise multiplication. According to Definition \ref{finitedual}
we have that $\A$ is a $k$-coalgebra with the following structure maps
$\Delta':\A\rightarrow \A\otimes\A$ such that $\Delta'(f)(x\otimes y) =
f(x\cdot y)$ for any $f\in\A$ and $\epsilon':\A\rightarrow k$ such that
$\epsilon'(f) = f(1)$, for any $f\in\A$.

\begin{corollary}
\label{grouplikebialgebra}
Let $(B,\Delta,\epsilon)$ be a $k$-bialgebra, then $k G(B^\circ)$ inherits
the structure of a $k$-bialgebra.
\end{corollary}
\proof
It is sufficient to note that $G(B^\circ)\subseteq B^\circ$ and hence the subspace
generated by $G(B^\circ)$, that is $k G(B^\circ)$, lies also in $B^\circ$. Indeed
$k G(B^\circ)$ is a subalgebra of $B^\circ$. If we denote by
$\Delta''$ the restriction of $\Delta'$ at $k G(B^\circ)$ and by $\epsilon''$
the restriction of $\epsilon'$ at $k G(B^\circ)$ as in the comments to Theorem
\ref{bialgebra}, we obtain that $(kG(B^\circ),\Delta'',\epsilon'')$ is a $k$-bialgebra. In fact, for any $f\in G(B^\circ)$, $\Delta''(f)=f\otimes f$ because $(C^\circ,\Delta',\epsilon')$ is a $k$-coalgebra. Then \[(\Delta''\otimes\id)(f\otimes f)=f\otimes f\otimes f=(\id\otimes\Delta'')(f\otimes f).\]Analogously, \[(\epsilon''\otimes\id)(f\otimes f)=1\otimes f=(\id\otimes\epsilon'')(f\otimes f)\]and we complete the proof.
\endproof

\begin{theorem}
\label{gradingaction}
Let $(S,\cdot)$ be a monoid and let $A$ be a $S$-graded algebra over $k$.
Let $H = k S$, then $A$ is a left $k G(H^\circ)$-module algebra.
\proof
In light of Corollary \ref{modulestructure}, one has that $A$ is a left $H^*$-module.
Indeed $k G(H^\circ)\subseteq H^*$ and therefore $A$ inherits the left
$k G(H^\circ)$-module structure. Corollary \ref{grouplikebialgebra} gives us
that $k G(H^\circ)$ is a $k$-bialgebra and hence it is enough to verify the conditions of a $k$-bialgebra being a module algebra. We observe that we may check only the case that $a,b\in A$ are homogenous elements
of degrees $s_a,s_b$ respectively. Then, for any $f\in G(H^\circ)$ we have that
$\Delta''(f) = f\otimes f$ by the definition of group like elements and
$\gamma(f\otimes a b) = f(s_a s_b)a b$. On the other hand,
$\gamma(f\otimes a)\gamma(f\otimes b) = f(s_a)f(s_b) ab$. By Proposition
\ref{grouplikealgebra} it turns out that $f$ is an algebra homomorphism, hence
$\gamma(f\otimes a)\gamma(f\otimes b) = \gamma(f\otimes ab)$. Moreover
$\gamma(f\otimes 1_A) = f(1_S)1_A = \epsilon''(f)1_A$ and we complete the proof.

\endproof
\end{theorem}

\begin{corollary}\label{classic1}
Let $A$ be a $S$-graded algebra over $k$, then $G((kS)^\circ)$ acts on $A$.
\end{corollary}
\proof
If $A$ is a $S$-graded algebra then $A$ is a left $kG((kS)^\circ)$-module algebra
by Theorem \ref{gradingaction}. The conclusion now follows from Proposition
\ref{actionmodulealgebra}.
\endproof

\section{Pontryagin duality}

From now on $k=\C$. In this section we introduce the essential topological tools in order to describe
a sub-bialgebra of $(k S)^\circ$ and its action on an $S$-graded algebra when $S$ is a bounded semilattice.

\begin{definition}
Let $(S,\cdot)$ be a monoid (group) endowed with the structure of a topological
space. We say that $S$ is a topological monoid (group) if the product $S\times S\to S$
is a continuous map, where $S\times S$ is endowed with the product topology.
\end{definition}

An important example of an abelian topological group is the compact 
set $\T =\{z\in\C \mid |z| = 1\}$. In the same way, one has the abelian compact
monoid $\B = \T \cup \{0\}$.

\begin{definition}
Let $k$ be a topological field and let $V$ be a finite dimensional $k$-vector
space endowed with the natural structure of topological space. Consider also
$\End_k(V)$ as a topological vector space and let $S$ be a topological
monoid. A continuous monoid homomorphism $\rho:S\rightarrow \End_k(V)$
is called a continuous linear (finite dimensional) representation of $S$. 
\end{definition}


Introducing topology one needs to modify slightly some algebraic structures such as group algebras (monoid algebras). In fact:
\begin{definition}
Let $S$ be a discrete monoid (group). The monoid (group) algebra $kS$ of $S$ is the vector space of all finitely supported $k$-valued functions on $S$, equipped with the convolution product $(f*g)(x):=\sum_{yz=x}f(y)g(z)$.
\end{definition}

\begin{remark}\label{solving}
We want to focus on the fact that group algebras of topological groups are defined by convolution products. In the case of discrete groups, the purely algebraic definition of group algebra and the previous one coincide. This is not the case of locally compact groups. In the last case the two definitions coincide if $G$ is finite.
\end{remark}

\begin{definition}
Let $(S,\cdot)$ be a topological monoid. Then, the set
\[
\widehat{S} =
\{\chi:S\rightarrow\B\mid \chi\ \mbox{continuous monoid homomorphism}\}
\]
is called the Pontryagin dual of $S$. Moreover, if $(G,\cdot)$ is a topological
group, then
\[
\widehat{G} =\{\chi:G\rightarrow\T\mid\chi\ \mbox{continuous monoid homomorphism}\}
\]
is called the Pontryagin dual of $G$.
\end{definition}

We remark that if $S$ is a monoid (group), then $\widehat{S}$ is a monoid (group) as well with the usual pointwise multiplication.

\begin{definition}
A topological monoid (group) is called discrete if it is a topological monoid
(group) endowed with the discrete topology.
\end{definition}

The following Duality Theorem is classical and can be found in \cite{mor1}:

\begin{theorem}
\label{discretecompact}
Let $G$ be a compact (discrete) topological group, then its Pontryagin dual
$\widehat{G}$ is discrete (compact) group.
\end{theorem}

We recall the following definition:

\begin{definition}
Let $(S,\cdot)$ be a topological monoid. We say that a continuous function $f:S\to k$
is representative if and only if the linear span of all left-translates of $f$
that are the functions
\[
f(z\ \cdot):y\mapsto f(zy)
\]
with $z\in S$, is finite dimensional. We shall denote the space of all
continuous representative functions on $S$ by $\rp(S)$.
\end{definition}

\begin{remark}\label{baobab}
By \cite{tim1}, Example 1.2.5, we have that the previous definition is equivalent to the following. A function $f:S\to k$ is representative if and only if there exists a
finite dimensional continuous representation $\rho:S\to\End_k(V)$ and elements
$v\in V$ and $\Phi\in V^*$ such that $f(s) = \Phi(\rho(s)v)$, for all $s\in S$.
\end{remark}

Now we denote by $k(S)$ the algebra of all $k$-valued functions on $S$,
where the addition and multiplication are defined pointwise. It is easy to note that $\rp(S)$ is a unital subalgebra of $k(S)$.
Define now the maps $\Delta:\rp(S)\rightarrow k(S\otimes S)$ and
$\epsilon:\rp(S)\rightarrow k$ as follows. For any $s,t\in S$ and $f\in\rp(S)$
we put $\Delta(f)(s\otimes t) = f(s\cdot t)$ and $\epsilon(f) = f(e)$. Using the same arguments as in \cite{tim1}, we have that $(\rp(S),\Delta,\epsilon)$ is a $k$-bialgebra.

Let $(S,\cdot)$ be a monoid. One denotes by $\rp_k(S)$ the set of representative
functions on $S$ without any hypothesis of continuity and over any field $k$.
Indeed, $\rp(S)\subseteq\rp_k(S)$. We consider now the following
result.

\begin{proposition}
If $S$ is a monoid and $H = k S$, then $H^\circ = \rp_k(S)$.
\end{proposition}

We would like to mention that the proof of the previous result is given for groups
in the book by Montgomery (\cite{mon1}), but it should be easily generalized to monoids
using Lemma 9.1.1 of \cite{mon1}.

An immediate consequence of Proposition \ref{comodulemodule} is the following.

%
%

We consider now the following result.

\begin{theorem}
\label{repdual}
Let $G$ be a discrete abelian group, then
\[
\rp(\widehat{G})\cong kG
\]
as $k$-bialgebras (in fact as Hopf algebras).
\end{theorem}

For a proof of the above result, see for instance the book of Timmermann \cite{tim1},
page 10. In the present paper, our main goal is to have a similar result for a special
class of monoids. For this purpose, let us consider the following definition.

\begin{definition}
Let $S$ be a monoid. We say that $S$ is an inverse monoid if for any $s\in S$ there exists
a unique $s^{-1}\in S$ such that $ss^{-1}s=s$ and $s^{-1}ss^{-1}=s^{-1}$.
\end{definition}

We recall that a set $S$ partially ordered by the binary relation $\leq$ is a \textit{semilattice} if for any $s,t\in S$, the greatest lower bound (glb) (or the smallest upper bound) of $\{s,t\}$ exists. If we set $s\leq t:=glb\{s,t\}$, it turns out that $(S,\leq)$ is a commutative nilpotent semigroup. We say that $(S,\leq)$ is a bounded semilattice if $(S,\leq)$ is a monoid. As an example, if $\leq$ denotes the usual comparing relation in $\Z$, then $(\Z,\leq)$ is a semilattice while $(\N,\leq)$ is a bounded semilattice.

Indeed, groups and bounded semilattices are examples of inverse monoids. The following results explain why the notion of inverse monoid is important (see \cite{bah1}
Lemma 8).

\begin{theorem}
\label{dualmonoid}
Let $S$ be a discrete abelian inverse monoid. Then $S\cong \widehat{\widehat{S}}$ as
topological monoids.
\end{theorem}


We have the following result.

\begin{theorem}
\label{repmonoid}
Let $S$ be a discrete bounded semilattice. Then, one has a bialgebras isomorphism
$\rp(\widehat{S})\cong k(S/I)$, for some biideal $I\subseteq kS$. In particular, if $S$
is finite, then $\rp(\widehat{S})\cong kS$ as bialgebras.
\end{theorem}

\proof
Indeed, one can identify $\widehat{S}$ with the monoid
\[
\{(z_s)_{s\in S}\in\prod_{s\in S}\B \mid z_s z_t = z_{st},\ \mbox{for all}\ s,t\in S\}
\subseteq \prod_{s\in S}\B.
\]
Since $\widehat{S}$ is a semilattice endowed with the pointwise multiplication, then
each of its continuous finite-dimensional representations is equivalent to a direct sum
of continuous one-dimensional ones. Such representations are in one-to-one correspondence
with elements $\chi$ of the Pontryagin dual of $\widehat{S}$, and hence, by Theorem
\ref{dualmonoid}, with elements of $S$.

By Remark \ref{baobab}, it turns out that the algebra $\rp(\widehat{S})$ is generated by the family of functions $(\ev_s)_{s\in S}$
given by $\ev_s(\chi) = \chi(s)$, for all $\chi\in\widehat{S}$. Now, we have that the map
$kS\rightarrow\rp(\widehat{S})$ such that $s\mapsto\ev_s$ for any $s\in S$ is a surjective
bialgebra homomorphism. In fact, one has
\[
\begin{array}{l}
(\ev_s\ev_t)(\chi) = \ev_s(\chi)\ev_t(\chi) = \chi(s)\chi(t) = \chi(st) = \ev_{st}(\chi); \\
\Delta(\ev_s)(\chi,\eta) = \ev_s(\chi\eta) = (\chi\eta)(s) = \chi(s)\eta(s) =
\ev_s(\chi)\ev_s(\eta); \\
\epsilon(\ev_s) = \ev_s(\widehat{e}) =\widehat{e}(s)=1.
\end{array}
\]
If $S$ is finite, the elements  of its Pontryagin
dual are exactly its characters in the usual sense. Hence, the functions
$(\ev_s)_{s\in S}$ are linearly independent (see for example Lemma III.2 of \cite{nee1})
and the map is a $k$-bialgebra isomorphism.
\endproof

\begin{corollary}
\label{classic2}
Let $S$ be a discrete bounded semilattice and let $A$ be a $S$-graded algebra.
Then, $\widehat{S}$ acts on $A$.
\end{corollary}
\proof
Due to the fact that $S$ is discrete we have $\rp_{k}(S) = \rp(S)$. Now Theorem
\ref{dualmonoid} says that $S\cong\widehat{\widehat{S}}$. Using Theorem \ref{repmonoid}
we have the following chain of isomorphisms:
\[
(kS)^\circ\cong\rp_k(S) = \rp(S_D)\cong\rp(\widehat{\widehat{S}}_D),\] where the lower index $D$ means that we are dealing with the discrete topology. We consider now $\widehat{S}_D$. By Theorem \ref{repmonoid}, we have $\rp(\widehat{\widehat{S}}_C)\cong k\widehat{S}/I$, for some biideal $I$ of $k\widehat{S}$, where we indicate by the lower index $C$ the compact topology (see Theorem \ref{discretecompact}). Indeed $\rp(\widehat{\widehat{S}}_C)\subseteq\rp(\widehat{\widehat{S}}_D)$, then we have:\[(kS)^\circ\cong\rp(\widehat{\widehat{S}}_D)\supseteq\rp(\widehat{\widehat{S}}_C)\cong k\widehat{S}/I.\]
In light of Corollary \ref{classic1} and Lemma \ref{grouplikequotient}, there exists an action $\rho$ of $\widehat{S}+I$ on $A$. Now the claim follows once we consider $\Gamma\circ\rho$, where $\Gamma$ is the canonical action of $\widehat{S}$ on $\widehat{S}+I$.
\endproof

\begin{remark}
If $G$ is a discrete finite abelian group, then $\rp_k(G) = \rp(G)$ and
$G\cong\widehat{\widehat{G}}$. Moreover, we have the following chain
of equalities:
\[
(kG)^\circ\cong\rp_k(G) = \rp(G)\cong\rp(\widehat{\widehat{G}})\cong k\widehat{G}.
\]
In light of Corollary \ref{classic1} and Remark \ref{solving}, if $A$ is a $G$-graded algebra, then the dual group $\widehat{G}$ acts on $A$.
\end{remark}

%

\section{Applications}

In this section we shall verify some of the results obtained in Sections 1
and 2 for concrete algebras.

\subsection{A finite monoid-grading on upper triangular matrices}
We consider the finite dimensional algebra $UT_2(k)$ of $2\times 2$ upper
triangular matrices with entries from $k=\C$. We consider the following two
subspaces of $UT_2(k)$:
\[
A_1=\left\{
\left(\begin{array}{cc}
	0 & 0\\
	0 & d
\end{array}\right)\left|\right. d\in k
\right\},\ \ 
A_2=\left\{
\left(\begin{array}{cc}
	a & b\\
	0 & 0
\end{array}\right)\left|\right. a,b\in k
\right\}.
\]
One has that $UT_2(k)=A_1\oplus A_2$. Let $X$ be a subset of $\N$ such that
$|X|=2$, say $X=\{n_1,n_2\}$, with $n_1<n_2$. If we set $\|A_1\| = n_1, \|A_2\| = n_2$,
we obtain that $UT_2(k)$ is a $(X,\max)$-graded algebra. In fact, one has:
\[
\begin{array}{l}
A_1\cdot A_2 =
\left(\begin{array}{cc}
	0 & 0\\
	0 & d
\end{array}\right)
\left(\begin{array}{cc}
	a & b\\
	0 & 0
\end{array}\right)=
\left(\begin{array}{cc}
	0 & 0\\
	0 & 0
\end{array}\right)\in A_2 =
A_{\max\{n_1,n_2\}},
\\ A_2\cdot A_1=\left(\begin{array}{cc}
	a & b\\
	0 & 0
\end{array}\right)
\left(\begin{array}{cc}
	0 & 0\\
	0 & d
\end{array}\right)=
\left(\begin{array}{cc}
	0 & bd\\
	0 & 0
\end{array}\right)\in A_2 =
A_{\max\{n_2,n_1\}},
\\ A_1\cdot A_1 =
\left(\begin{array}{cc}
	0 & 0\\
	0 & d
\end{array}\right)
\left(\begin{array}{cc}
	0 & 0\\
	0 & x
\end{array}\right)=
\left(\begin{array}{cc}
	0 & 0\\
	0 & dx
\end{array}\right)\in A_1 =
A_{\max\{n_1,n_1\}},
\\ A_2\cdot A_2=\left(\begin{array}{cc}
	a & b\\
	0 & 0
\end{array}\right)\left(\begin{array}{cc}
	x & y\\
	0 & 0
\end{array}\right)=\left(\begin{array}{cc}
	ax & ay\\
	0 & 0
\end{array}\right)\in A_2 =
A_{\max\{n_2,n_2\}}.
\end{array}
\]
Let $H = k X$ and consider the set $O =\{f_1,f_2\}$,
where $f_i:H\rightarrow k$ are homomorphisms of algebras such that
\[
f_1(x)=\left\{
\begin{array}{cc}
1 & \mbox{if}\ x=n_1;\\
0 & \mbox{if}\ x=n_2;
\end{array}
\right.\ \ \
f_2(x)=\left\{
\begin{array}{cc}
1 & \mbox{if}\ x=n_1;\\
1 & \mbox{if}\ x=n_2.
\end{array}
\right.
\]
If $f\in H^*$, then $f\in\alg$ if $f(n_2)=f(\max\{n_1,n_2\})=f(n_1)f(n_2)$. If $f(n_2)\neq0$,
then $f(n_1) = 1$. In this case $f(n_2)=f(n_2)^2$ that is $f(n_2)=1$ and we conclude
that $f=f_i$ for $i\in\{1,2\}$. Moreover, one has
$(\alg,\cdot)=(O,\cdot)\cong(X,\min)$. By Corollary
\ref{classic1} we have that there exists an action of $(X,\min)$ on $UT_2(k)$
such that if
$A=\left(\begin{array}{cc}
	a & b\\
	0 & d
\end{array}\right),$ then
\[
\gamma\left(n_1\otimes A\right) =
\left(\begin{array}{cc}
	0 & 0\\
	0 & d
\end{array}\right)
\ \mbox{and}\ 
\gamma\left(n_2\otimes A\right) = A.
\]
Note that the previous grading may be extended to the algebra $UT_m(k)$ of upper
triangular matrices of size $m$ in a natural way. More precisely, let $X$ be a
subset of $\N$ such that $|X| = m$, say $X=\{n_1,n_2,\ldots,n_m\}$ with
$n_1<n_2<\cdots<n_m$. Again $(X,\max)$ is an abelian semilattice
and, if we define
\[
A_1=\left\{\left(\begin{array}{ccccc}
	0 & \ldots & \ldots & 0\\
	\vdots & & &\vdots\\
	0 & \ldots & \ldots & 0\\
	0 & \ldots & 0 & a
\end{array}\right)\right\},
A_2=\left\{\left(\begin{array}{ccccc}
	0 & \ldots & \ldots & 0\\
	\vdots & & &\vdots\\
	0 & \ldots & b & c\\
	0 & \ldots & \ldots & 0
\end{array}\right)\right\},
\ldots,
\] 
where $\|A_i\| = n_i$, we obtain that $UT_m(k) = \bigoplus_{i=1}^m A_i$ is a
$(X,\max)$-graded algebra. Then, one obtains similar results as in the case of $UT_2(k)$. 

If we endow $X$ with the discrete topology, then $(X,\max)$ is a discrete bounded semilattice. Indeed $(X,\min)\cong\widehat{(X,\max)}\cong(X,\max)$ as we have shown in Corollary \ref{classic2}.

\subsection{The case of algebra of letterplaces with a $(\nnn,\max)$-grading}
We apply the results of Sections 1 and 2 to the algebra of letterplaces.
We mention that a detailed exposition of this subject can be found in the monograph
by Grosshans, Rota and Stein \cite{grs1}.

We start with the following definition. In what follows, we assume that the base field $k$
is of characteristic zero. Let $X = \{x_1,x_2,\ldots\}$ be a finite or countable set
endowed with a map $|\cdot|:X\to \Z_2$. Put $X_k = \{x_i\in X\mid |x_i| = k\}$ for $k=0,1$
and hence $X = X_0\cup X_1$. We denote by $F = k\langle X\rangle$ the free associative
algebra generated by $X$ that is the tensor algebra of the vector space $k X$. This algebra
is clearly $\Z_2$-graded if we put $|w| = |x_{i_1}| + \cdots + |x_{i_n}|$, for any
monomial $w = x_{i_1} \cdots x_{i_n}$. Note that a $\Z_2$-graded algebra is called
also a superalgebra. Now, denote by $I$ the two-sided ideal of $F$ generated the binomials
\[
x_i x_j - (-1)^{|x_i||x_j|} x_j x_i,
\]
where $x_i,x_j\in X$. Clearly $I$ is $\Z_2$-graded ideal and we define $\Super[X] = F/I$.
This $\Z_2$-graded algebra is called the free supercommutative algebra generated by $X$.
The identities of $\Super[X]$ defined by the above binomials are said in fact supercommutative
identities. It is plain that a $k$-basis of $\Super[X]$ is given by the cosets of the monomials
$w = x_{i_1} \ldots x_{i_n}$ such that
\[
x_{i_k}\leq x_{i_{k+1}},\ \mbox{with}\ x_{i_k} = x_{i_{k+1}}\ \mbox{only if}\
|x_{i_k}| = 0.
\]
Note that the algebra $\Super[X]$ is isomorphic to the tensor product
$k[X_0] \otimes E[X_1]$ where $k[X_0]$ is the polynomial algebra in the commuting
variables $x_i\in X_0$ and $E[X_1]$ is the exterior (or Grassmann) algebra of the
vector space $k X_1$.

Consider now set of positive integers $\N^*$ be endowed with a map $|\cdot|:\N^*\to\Z_2$.
Then, consider the product set $X\times\N^*$ and define the map
$|\cdot|:X\times \N^*\to\Z_2$ such that $|(x_i,j)| = |x_i| + |j|$, for any $i,j\geq 1$.
We put $\Super[X|\N^*] = \Super[X \times \N^*]$ and we call this supercommutative algebra
the letterplace superalgebra. A variable $(x_i,j)\in X \times \N^*$ will be denoted as
$(x_i|j)$ where $x_i$ is said  a ``letter'' and $j$  a ``place''. These names arise
from the possibility to embed the tensor superalgebra $k\langle X\rangle$ into
$\Super[X \times \N^*]$ by using the map $w = x_{i_1} \ldots x_{i_n}\mapsto
m = (x_{i_1}|1)\ldots (x_{i_n}|n)$.

Define $\overline{\N} = \N\cup\{-\infty\}$ and consider the discrete abelian monoid
$(\overline{\N},\max)$ with neutral element $-\infty$. Note that $(\overline{\N},\max)$
is in fact a semilattice. Then, one can introduce the following grading of $P =
\Super[L \times \N]$ over the monoid $(\overline{\N},\max)$. Denote by $M$ the set
of monomials of $P$.

\begin{definition}
Let $\w:M\to\nnn$ be the unique mapping such that
\begin{itemize}
\item[(i)] $\w(1) = -\infty$;
\item[(ii)] $\w(m n) = \max(\w(m),\w(n))$, for any $m,n\in M, m n\neq 0$;
\item[(iii)] $\w((x_i|j))) = j$, for all $i,j\geq 1$.
\end{itemize}
We call $\w$ the weight function of $P$. If $P_{(i)}\subset P$ is the subspace spanned
by all monomials of weight $i$ then $P = \bigoplus_{i\in\nnn} P_{(i)}$ is grading of $P$
over the semilattice $(\nnn,\max)$.
\end{definition}

Note that the idea of such a grading when $P$ is a commutative polynomial algebra
has been introduced in \cite{LSL}. We provide $H = k\overline{\N}$ with the usual structure of a $k$-bialgebra. From Theorem
\ref{bialgebra} we have that $H^\circ$ is also a $k$-bialgebra. Proposition
\ref{comodulegrading} gives us that $P$ is a right $H$-comodule under the map
$\rho:P\rightarrow P\otimes H$ such that $\rho(p) = \sum_{n\in\overline{\N}}p_n\otimes n$,
for any $p = \sum_{n\in\overline{\N}}p_n\in P$. Furthermore, Proposition
\ref{comodulemodule} gives us that $P$ is a left $H^*$-module under the action
$\gamma:H^*\otimes P\rightarrow P$ such that $\gamma(f\otimes p) =
\sum_{n\in\overline{\N}}f(n)p_n$, for any $f\in H^*$. Due to the fact that
$\alg = G(H^\circ)\subseteq H^\circ\subseteq H^*$, we obtain the action of
$G(H^\circ)$ over $P$. Let $\underline{\N} = \N\cup\{+\infty\}\cup\{-\infty\}$
and define the set $O = \{f_n|n\in\underline{\N}\}$ as follows. For any
$n\in\overline{\N}$, the function $f_n:H\rightarrow k$ is such that for all
$m\in\overline{\N}$
\[
f_n(m)=\left\{
\begin{array}{ll}
0 & \mbox{if}\ m > n; \\
1 & \mbox{otherwise}.
\end{array}
\right.
\]
Moreover, we define $f_{+\infty}:H\rightarrow k$ such that $f_{+\infty}(m) = 1$,
for all $m\in\overline{\N}$.

\begin{proposition}\label{appl}
For any $n\in\underline{\N}$, $f_n\in\alg$.
\end{proposition}
\proof
For any $a,b\in \overline{\N}$, consider $f_n(\max\{a,b\})$. We distinguish
three cases. If $a,b>n$ and $a\geq b$, then $f_n(\max\{a,b\}) = f_n(a) = 0$
and $f_n(a)f_n(b)=0$. Moreover, if $a>n\geq b$, then $f_n(\max\{a,b\})=f_n(a)=0$
and $f_n(a)f_n(b)=0$. Finally, if $b\leq a\leq n$, then $f_n(\max\{a,b\})=f_n(a)=1$
and $f_n(a)f_n(b)=1$.
\endproof

\begin{proposition}
$O=\alg$, with $(O,\cdot)\cong(\underline{\N},\min)$.
\end{proposition}
\proof
Proposition \ref{appl} gives us that $O\subseteq\alg$. We are going to prove the reverse
inclusion. If $f\in H^*$, then $f\in\alg$ if and only if for any $a,b\in\nnn$
such that $a\geq b$ one has $f(a)=f(\max\{a,b\})=f(a)f(b)$. If $f(a)\neq0$,
then $f(b) = 1$. Moreover, in this case $f(a)=f(a)^2$ that is $f(a)=1$ and we conclude
that $f=f_c$ for some $c\geq a$. We consider now the monoid $(O,\cdot)$ endowed
with the pointwise multiplication of functions. Let $s,t\in\underline{\N}$ such that
$s\leq t$. Then $f_s\cdot f_t = f_{\min\{s,t\}}$. In fact, for any $n\in\N$ one has
\[
(f_s\cdot f_t)(n)=\left\{
\begin{array}{cc}
0 & n>t\geq s \\
0 & t\geq n>s \\ 
1 & t\geq s\geq n
\end{array}
\right.
\]
and
\[
f_s(n)=\left\{
\begin{array}{cc}
0 & n>s \\
1 & n\leq s.
\end{array}
\right.
\]
Note that $(O,\cdot)$ is a monoid with neutral element $f_{+\infty}$. Then, the
correspondence $z\leftrightarrow f_z$ gives us $(O,\cdot)\cong(\underline{\N},\min)$
as monoids.
\endproof
 
We shall investigate now the structure of $H^\circ$. We start with the following lemma.

\begin{lemma}
\label{matrix}
Let $n\geq2$ and consider the following $n\times n$ matrix with entries in $\C$
\[
M =\left(
\begin{array}{ccccc}
a_{11} & a_{12} & a_{13} &\ldots & a_{1n}\\
a_{12} & a_{12} & a_{13} &\ldots & a_{1n}\\
\vdots & & & &\vdots\\
a_{1n} & a_{1n} & a_{1n} &\ldots & a_{1n}\\
\end{array}
\right),
\]
where $a_{1n}\neq0$ and for any $i=1,\ldots,n$ one has $a_{1i}\neq a_{1i+1}$.
Then $\det(M)\neq 0$.
\end{lemma}
\proof
We proceed by induction on $n$. If $n=2$, then
$M=
\left(
\begin{array}{cc}
a_{11} & a_{12}\\
a_{12} & a_{12}
\end{array}
\right)$ is a $2\times2$ matrix such that $a_{12}\neq0$ and $a_{11}\neq a_{12}$. Then
$\det(M)=a_{11}a_{12}-a_{12}^2\neq0$ and the conclusion follows.
Consider now the $(n+1)\times (n+1)$ matrix
\[
M =\left(
\begin{array}{ccccc}
a_{11} & a_{12} & a_{13} &\ldots & a_{1n+1}\\
a_{12} & a_{12} & a_{13} &\ldots & a_{1n+1}\\
\vdots & & & &\vdots\\
a_{1n+1} & a_{1n+1} & a_{1n+1} &\ldots & a_{1n+1}\\
\end{array}
\right),
\]
where $a_{1n+1}\neq 0$ and $a_{1i}\neq a_{1i+1}$ for all $i=1,\ldots,n+1$.
We observe that $\det(M) = \det(M')$, where 
\[
M'=\left(
\begin{array}{ccccc}
a_{11} & a_{12} & a_{13} &\ldots & a_{1n+1}\\
a_{12}-a_{11} & 0 & 0 &\ldots & 0\\
\vdots & & & &\vdots\\
a_{1n+1} & a_{1n+1} & a_{1n+1} &\ldots & a_{1n+1}\\
\end{array}
\right).\]It turns out that $\det(M')=(a_{11}-a_{12})\det(M''),$ where \[M''=\left(
\begin{array}{ccccc}
a_{12} & a_{13} & a_{14} &\ldots & a_{1n+1}\\
a_{13} & a_{13} & a_{14} &\ldots & a_{1n+1}\\
\vdots & & & &\vdots\\
a_{1n+1} & a_{1n+1} & a_{1n+1} &\ldots & a_{1n+1}\\
\end{array}
\right).\]
Now the inductive step applies at $M''$ that is $\det(M'')\neq 0$ and the conclusion
follows because $a_{11}\neq a_{12}$. 
\endproof

\begin{proposition}
\label{dual2}
Let $f\in H^*$. Then $f\in H^\circ$ if and only if $f$ assumes definitely constant
values on $\nnn$.
\end{proposition}
\proof
If $f\in H^\circ$ then $\dim_k\{\nnn\cdot f\}<\infty$. Hence there exist
$s_1,\ldots,s_l\in\nnn$ such that, for any $n\in\nnn$, there exist
$\alpha_1^n,\ldots,\alpha_l^n$ such that
\[
n\cdot f=\sum_{i=1}^l\alpha_i^n(s_i\cdot f).
\]
Then, for any $m\in\nnn$ one has
\[
(n\cdot f)(m)=\sum_{i=1}^l\alpha_i^nf(\max\{s_i,m\}).
\]
Consider the case $m=n$, where $n>s_l>\cdots>s_1$. We have
\[
f(n)=\left(\sum_{i=1}^l\alpha_i^n\right)f(n).
\]
If $\sum_{i=1}^l\alpha_i^n=0$, then $f(n)=0$, otherwise $\sum_{i=1}^l\alpha_i^n=1$.
If $m>n$, then
\[
f(n)=\left(\sum_{i=1}^l\alpha_i^n\right)f(m)
\]
that is $f(n)=f(m)$. The latter result implies that if $f$ does not assume definitely
the value 0 on $\nnn$ and there exist $m_1,m_2\in\nnn$ such that $f(m_1)\neq 0$ and
$f(m_2)\neq 0$. Therefore, if $m_1<m_2$ one has that $f(m_1) = f(m_2)$ and the proof
follows. Conversely, if $f$ assumes definitely constant values on $\nnn$
we have $f(\nnn)=\{a_1,\ldots,a_l\}$, where $a_1<a_2<\cdots<a_l$. Let us consider
the following element of $\nnn$:
\[
s_1=\max\{p\mid f(p)=a_1,f(r)=a_1\ \mbox{if}\ r\leq p, f(p+1)\neq a_1\}.
\]
Then, we consider
\[
s_2=\max\{p\mid f(p)=a_i,f(r)=a_i\ \mbox{if}\ s_1<r\leq p, f(p+1)\neq a_i\}
\]
and so on. Note that these elements are finite in number, then we define
\[
S=\{s_1,s_2,\ldots,s_t\},
\]
where $s_1 < s_2 <\cdots < s_t$. We claim that $\{S\cdot f\}$ is a basis of $\nnn\cdot f$.
First of all, we show that $S$ is a linearly independent set. Suppose that there
exist $\alpha_1,\ldots,\alpha_t\in k$ such that $\sum_{i=1}^t\alpha_i(s_i\cdot f) = 0$.
Then, for any $m\in\nnn$ we have
\[
\sum_{i=1}^t\alpha_if(\max\{s_i,m\})=0.
\]
If we substitute $m$ with each of the $i$'s, we obtain the following linear system
in the indeterminates $\alpha_1,\ldots,\alpha_t$
\[
\left\{\begin{array}{l}
	\alpha_1a_1+\alpha_2a_2+\cdots+\alpha_ta_t=0\\
	\alpha_1a_2+\alpha_2a_2+\cdots+\alpha_ta_t=0\\
	\vdots\\
	(\alpha_1+\cdots+\alpha_t)a_t=0.
\end{array}\right.
\]
We may always suppose that $a_t\neq 0$. Then, by Lemma \ref{matrix}, the previous system
has only the trivial solution and the claim follows. For the generation property we have
to prove that for any $n\in\nnn$ there exist $\alpha_1^n,\ldots,\alpha_l^n$ such that
for any $m\in\nnn$ one has that
\[
f(\max\{n,m\})=\sum_{i=1}^l\alpha_i^n(f(\max\{s_i,m\})).
\]
For simplicity, we consider only the case $n>s_t>\cdots>s_1$ because the remaining cases
can be treated similarly. If we substitute $m$ with each of the $i$'s, we obtain the
following linear system in the indeterminates $\alpha_1,\ldots,\alpha_t$:
\[
\left\{\begin{array}{l}
	\alpha_1a_1+\alpha_2a_2+\cdots+\alpha_ta_t=a_t\\
	\alpha_1a_2+\alpha_2a_2+\cdots+\alpha_ta_t=a_t\\
	\vdots\\
	(\alpha_1+\cdots+\alpha_t)a_t=a_t.
\end{array}\right.\]
Now the conclusion follows from Lemma \ref{matrix}. 
\endproof

The next theorem relates $H^\circ$ with its group-like elements.

\begin{theorem}
If $f\in H^\circ$ then $f$ is a linear combination of elements of $G(H^\circ)$.
\end{theorem}
\proof
By Proposition \ref{dual2}, we have that $f$ assumes only a finite number of values on $\nnn$,
that is $\{a_1,\ldots,a_l\}$. We set
\[
P=|\{s\mid f(s)=a_1,f(r)=a_1\ \mbox{if}\ r<s\ \mbox{and}\ f(s+1)\neq a_1\}|.
\]
Now we identify $f$ with the vector
\[
v=(\underbrace{a_1,\ldots,a_1}_{P},a_2,\ldots,a_{l-1},a_l,\ldots,a_l,\ldots)
\]
of the values that $f$ assumes on $\nnn$. Indeed, every $f_n\in G(H^\circ)$ may be identified
with the vector
\[
v_n=(\underbrace{1,\ldots,1}_{P_n},0,\ldots,0,\ldots),
\]
where
\[
P_n=|\{s\mid f_n(s)=1,f_n(r)=1\ \mbox{if}\ r<s\ \mbox{and}\ f(s+1)=0\}|.
\]
Note that $f_{+\infty}$ may be identified with the vector
\[
v_{+\infty}=(1,\ldots,1,\ldots).
\]
Simple linear algebra considerations show that the set
\[
\{v_{{+\infty}},v_{P+l-1},v_{P+l-2},\ldots,v_{P+1}\}
\]
is a linearly independent generating set and the claim follows.
\endproof

\begin{corollary}\label{final}
$H^\circ\cong k\otimes kG(H^\circ)\cong kG(H^\circ).$
\end{corollary}

Finally, there exists an action of $(\underline{\N},\min)$ on $P$ such that for any
$z\in \underline{\N}$
\[
\gamma(z\otimes p)=\sum_{n\in\overline{\N}}f_z(n)p_n.
\]
In particular, if $(x_i|j)$ is any variable of $P$ then it has degree $j$ and
then
\[
\gamma(z\otimes (x_i|j))=\left\{
\begin{array}{cc}
0 & j > z \\
(x_i|j) & j\leq z.
\end{array}\right.
\]

If we endow $\overline{\N}$ with the discrete topology, then $(\overline{\N},\max)$ is a discrete bounded semilattice. Indeed $(\underline{\N},\min)\cong\widehat{(\overline{\N},\max)}$ as we have shown in Corollary \ref{classic2} but, in this case, $(\overline{\N},\max)$ and $(\underline{\N},\min)$ are not isomorphic monoids.

\begin{remark}
In the proof of Theorem \ref{repmonoid} the biideal $I$ is the zero biideal
if and only if the family of functions $\{\ev_s\}_{s\in S}$ is linearly independent.
Following the example at page 10 of \cite{tim1}, the author proves the linearly independence
of the family $(\ev_g)_{g\in G}$ using a specific Haar measure on $\widehat{G}$.
In the literature we have not a general method to build up Haar measures on a monoid
structure. 
\end{remark}

%
%

\end{document}